\documentclass{amsart}
\usepackage{amssymb}
\usepackage[all]{xy}
\newcommand\datver[1]{\def\datverp%
  {\par\boxed{\boxed{\text{#1; Run: \today}}}}}

\newcommand\boxb[1]{\square_b}
\newcommand\Ahat{\widehat{A}}

\newcommand\paperbody%
         {}

\newtheorem{lemma}{Lemma}
\newtheorem{proposition}{Proposition}

\newtheorem{theorem}{Theorem}
\newtheorem{non-theorem}{Non-Theorem}

\theoremstyle{remark}
\newtheorem{definition}{Definition}
\newtheorem{remark}{Remark}

\newcommand\coF{{}^{\mathcal{C}}\kern-2pt\Lambda}

\newcommand\cFTs{{}^{\Phi}\overline{T}\kern-1pt{}^*}

\newcommand\even{\text{even}}
\newcommand\odd{\text{odd}}

\newcommand\Tr{\operatorname{Tr}}

\newcommand\bTr{\overline{\operatorname{Tr}}}
\newcommand\STr{\operatorname{STr}}

\newcommand\rTr{\operatorname{Tr_{R}}}

\newcommand\Cl{\operatorname{Cl}}
\newcommand\CCl{\CC\text{l}}
\newcommand\eCCl{\widehat{\CC\text{l}}}

\newcommand\SO{\operatorname{SO}}

\newcommand\Ch{\operatorname{Ch}}
\newcommand\Td{\operatorname{Td}}
\newcommand\dbar{\overline{\pa}}

\newcommand\hL{\widehat{L}}
\newcommand\tA{\tilde{A}}

\newcommand\tE{\tilde{E}}
\newcommand\tF{\tilde{F}}
\newcommand\tG{\tilde{G}}
\newcommand\tL{\tilde{L}}
\newcommand\tS{\tilde{S}}

\newcommand\tV{\tilde{V}}

\newcommand\com[1]{\overline{#1}}

\newcommand\cf{cf\@. }
\newcommand\ie{i\@.e\@. }

\newcommand\ecA{\widehat{\mathcal{A}}}

\newcommand\cA{\mathcal{A}}

\newcommand\cF{\mathcal{F}}
\newcommand\cH{\mathcal{H}}
\newcommand\cK{\mathcal{K}}

\newcommand\cL{\mathcal{L}}
\newcommand\cP{\mathcal{P}}
\newcommand\cQ{\mathcal{Q}}
\newcommand\cW{\mathcal{W}}

\newcommand\CC{\mathbb C}

\newcommand\ZZ{\mathbb Z}

\newcommand\bbC{\mathbb C}

\newcommand\bbQ{\mathbb Q}
\newcommand\bbR{\mathbb R}

\newcommand\bbZ{\mathbb Z}

\newcommand\cC{\mathcal C}

\newcommand\CIc{{\mathcal{C}}^{\infty}_{\text{c}}}
\newcommand\Kc{K_{\text{c}}}

\newcommand\CI{{\mathcal{C}}^{\infty}}

\newcommand\CmI{{\mathcal{C}}^{-\infty}}

\newcommand\Diag{\operatorname{Diag}}

\newcommand\Diff[1]{\operatorname{Diff}^{#1}}

\newcommand\cFNs{{}^{\Phi}\overline N\kern-1pt{}^*}

\newcommand\ind{\operatorname{ind}}
\newcommand\inda{\ind_{\text{a}}}

\newcommand\tr{\operatorname{tr}}

\newcommand\Hom{\operatorname{Hom}}
\newcommand\Id{\operatorname{Id}}

\newcommand\SU{\operatorname{SU}}
\newcommand\PU{\operatorname{PU}}
\newcommand\UU{\operatorname{U}}

\newcommand\dCI{\dot{\mathcal{C}}^{\infty}}

\newcommand\pa{\partial}

\newcommand\Spin{\operatorname{Spin}}
\newcommand\SpinC{\operatorname{Spin}_{\bbC}}
\newcommand\spin{\operatorname{spin}}
\newcommand\spinC{\operatorname{spin}_{\bbC}}

\newcommand\cl{\operatorname{cl}}

\renewcommand\Re{\operatorname{Re}}

\newcommand\Mif{\text{ if }}

\newcommand\Mor{\text{ or }}

\datver{1.2K; Revised: 24-8-2005}
\begin{document}
\title[Fractional analytic index]
{Fractional analytic index}

\author{V. Mathai}
\address{Department of Mathematics, University of Adelaide,
Adelaide 5005, Australia}
\email{vmathai@maths.adelaide.edu.au}
\author{R.B. Melrose}
\address{Department of Mathematics,
Massachusetts Institute of Technology,
Cambridge, Mass 02139}
\email{rbm@math.mit.edu}
\author{I.M. Singer}
\address{Department of Mathematics,
Massachusetts Institute of Technology,
Cambridge, Mass 02139}
\email{ims@math.mit.edu}
\begin{abstract}
For a finite rank projective bundle over a compact manifold, so associated
to a torsion, Dixmier-Douady, 3-class, $w,$ on the manifold, we define the
ring of differential operators `acting on sections of the projective bundle' in a
formal sense. In particular, any oriented even-dimensional manifold carries
a projective spin Dirac operator in this sense. More generally the
corresponding space of pseudodifferential operators is defined, with
supports sufficiently close to the diagonal, \ie the identity relation. For
such elliptic operators we define the numerical index in an essentially
analytic way, as the trace of the commutator of the operator and a
parametrix and show that this is homotopy invariant. Using the 
heat kernel method for the twisted, projective spin Dirac operator,
we show that this index is given by the usual formula, now in terms of the
twisted Chern character of the symbol, which in this case defines an
element of K-theory twisted by $w;$ hence the index is a rational number but in
general it is not an integer.
\end{abstract}
\maketitle


The Atiyah-Singer index theorem for an elliptic (pseudodifferential)
operator gives an integrality theorem; namely a certain characteristic
integral is an integer because it is the index of an elliptic
operator. Notably, for a closed spin manifold $Z$, the $\Ahat$ genus,
$\int_Z {\widehat A}(Z)$ is an integer because it is equal to the index of
the Dirac operator on $Z.$ When $Z$ is not a spin manifold, the spin bundle
$S$ does not exist, as a vector bundle, and when $Z$ has no $\spinC$
structure, there is no global vector bundle resulting from the patching of
the local bundles $S\otimes\mathcal L_i,$ where the $\mathcal L_i$ are line
bundles. However, as we show below, $S$ is a always a projective vector
bundle associated to the the Clifford algebra $\CCl(Z),$ of the cotangent
bundle $T^*Z$ which is an Azumaya bundle \cf
\cite{Mathai-Melrose-Singer1}. Such a (finite rank) projective vector
bundle, $E,$ over a compact manifold has local trivializations which may
fail to satisfy the cocycle condition on triple overlaps by a scalar
factor; this defines the Dixmier-Douady invariant in $H^3(Z,\bbZ).$ If this
torsion twisting is non-trivial there is no, locally spanning, space of
global sections. The Dixmier-Douady invariant for $\CCl(Z)$ is the third
integral Stieffel-Whitney class, $W_3(Z).$ In particular, the spin Dirac
operator does not exist when $Z$ is not a spin manifold. Correspondingly
the $\Ahat$ genus is a rational number, but not necessarily an integer. In
this paper we show that, in the oriented even-dimensional case, one can
nevertheless define a projective spin Dirac operator, with an analytic
index valued in the rational numbers, and prove the analogue of the
Atiyah-Singer index theorem for this operator twisted by a general
projective bundle. In fact we establish the analogue of the Atiyah-Singer
index theorem for a general projective elliptic pseudodifferential
operator. In a subsequent paper the families case will be discussed.

For a compact manifold, $Z,$ and vector bundles $E$ and $F$ over $Z$ the
Schwartz kernel theorem gives a one-to-one correspondence between
continuous linear operators from $\cC^\infty(Z, E) $ to $\CmI(Z,F)$ and
distributions in $\CmI(Z^2,\Hom(E, F)\otimes\Omega_R).$ Here $\Hom(E, F)$
is the `big' homomorphism bundle over $Z^2$ with fibre at $(z,z')$ equal to
$\hom(E_{z'},F_z)\equiv F_z\otimes E_{z'}',$ and $\Omega_R$ is the density
bundle lifted from the right factor. Restricted to pseudodifferential
operators of order $m,$ this becomes an isomorphism to the space,
$I^m(Z^2,\Diag;\Hom(E,F)\otimes\Omega _R),$ of conormal distributions with
respect to the diagonal, \cf \cite{Hormander3}.

This fact motivates our definition of projective pseudodifferential
operators when $E$ and $F$ are only projective vector bundles associated to
a fixed finite-dimensional Azumaya bundle $\cA.$ The homomorphism bundle
$\Hom(E,F)$ is then again a projective bundle on $Z^2$ associated to the
tensor product $\cA_L\otimes\cA_R'$ of the pull-back of $\cA$ from the left
and the conjugate bundle from the right. In particular if $E$ and $F$ have
DD invariant $\tau\in H^3(Z;\bbZ)$ then $\Hom(E,F)$ has DD invariant
$\pi_L^*\tau-\pi_R^*\tau\in H^3(Z^2;\bbZ).$ Since this class is trivial in
a tubular neighborhood of the diagonal it is reasonable to expect that
$\Hom(E,F)$ may be realized as an ordinary vector bundle there. In fact
this is the case and there is a canonical choice, $\Hom^{\cA}(E,F)$ of
extension. This allows us to identify the space of projective 
pseudodifferential operators, with kernels supported in a sufficiently
small neighborhood $N_\epsilon$ of the diagonal, 
$\Psi^\bullet_\epsilon(Z; E, F)$ with the space of conormal distributions
$I^\bullet_\epsilon (N_\epsilon,\Diag;\Hom^{\cA}(E, F)\otimes\Omega _R).$
Despite not being a space of operators, this has precisely the same local
structure as in the untwisted case and has similar composition properties
provided supports are restricted to appropriate neighborhoods of the
diagonal. The space of projective smoothing operators
$\Psi^{-\infty}_\epsilon(Z; E, F)$ is therefore identified with
$\CIc(N_\epsilon;\Hom^{\cA}(E, F)\otimes \pi_R^*\Omega).$ The principal
symbol map is well defined for conormal distributions so this leads
directly to the symbol map on $\Psi^m_\epsilon(Z; E, F)$ with values in
smooth homogeneous sections of degree $m$ of $\hom(E,F),$ the `little' or
`diagonal' homomorphism bundle which \emph{is} a vector bundle. Thus
ellipticity is well defined, as the invertibility of this symbol. The
`full' symbol map is given by the map to the quotient
$\Psi^\bullet_\epsilon(Z; E, F)/\Psi^{-\infty}_\epsilon(Z; E, F).$ The
usual calculus can then be applied and ellipticity, as invertibilty of the
principal symbol, implies invertibility of the image of an operator in this
quotient. Any lift, $B,$ of the inverse is a parametrix for the given
elliptic operator, $A.$ The analytic index of the projective elliptic
operator is then defined by
\begin{equation}
\inda(A)=\Tr(AB-\Id_F)-\Tr(BA-\Id_E)=\Tr([A,B])
\label{mms3.49a}\end{equation}
where the last expression, though compact, is slightly misleading. Directly
from this formula it appears that $\inda(A)$ might be complex valued. Using
the homotopy invariance discussed below it can be see directly to be real;
from the index formula it follows that $\inda(A)$ is rational.

In fact we further extend this discussion by allowing twisting by line
bundles defined over the bundle of trivializations of the Azumaya bundle;
these are $N$th roots of line bundles over the base and have actions of
$\SU(N)$ arising from fibre-holonomy. This allows us to include in the same
general framework the case of $\spinC$ Dirac operators on non-spin
manifolds. Thus we define projective $\spinC$ Dirac operators even when
there is no $\spinC$ structure.

Within the projective pseudodifferential operators, acting between two
projective bundles associated to the same Azumaya bundle, there is always a
full algebra of differential operators, with kernels supported within the
diagonal. On an even-dimensional oriented manifold the Clifford bundle is
an Azumaya bundle and has associated to it projective spin bundles, $S^\pm.$
The choice of a compatible connection gives a projective spin Dirac
operator. Such a projective Dirac operator, $\eth,$ can be coupled to
any unitary projective vector bundle $E$ over $Z$ associated to an Azumaya
algebra $\cA.$ Thus $S\otimes E$ is a projective vector bundle associated
to the Azumaya algebra $\CCl(Z)\otimes\cA.$ This coupled operator,
$\eth^+_E,$ is elliptic and its analytic index, in the sense defined above,
is
$$
\int_Z {\widehat A}(Z)\wedge\Ch_\cA(E)
$$
where $\Ch_\cA : K^0(Z, \cA) \to H^{even}(Z, \mathbb Q)$ is the twisted
Chern character.

When $Z$ is even dimensional, $\Kc^0(T^*Z,\pi^*\cA) \otimes \mathbb Q$ is
generated by the classes of symbols of such coupled signature
operators. We conclude from this, essentially as in the untwisted case,
that for a general projective elliptic pseudodifferential operators $T\in
\Psi^\bullet_\epsilon(Z; E, F)$ with principal symbol $\sigma_m(T),$
$$
\inda(T) =\int_{T^*Z}\Td(T^*Z)\wedge\Ch_{\cA}(\sigma_m(T)).
$$
This in turn shows the rationality of the analytic index and we conclude by
providing several examples where $\ind_a$ is not an integer, but only a
fraction, as justification for the title of the paper.

In the first section below the `big' homomorphism bundle is constructed,
near the diagonal, for any two projective bundles associated to the same
Azumaya bundle. In Section~\ref{P-twisting} we discuss a more general
construction of such `big' homomorphism bundles which corresponds to
twisting by a line bundle, only the $N$th power of which is well-defined
over the base. The projective spin bundle is discussed in the next section
as is its relationship to the spin bundle when a spin structure exists and
to $\spinC$ bundles when they exist. In Sections~\ref{Smooth} and
\ref{Pseudo} the notion of smoothing, and then general pseudodifferential,
operators between projective bundles (and more generally twisted projective
bundles) is introduced and for elliptic operators the index is defined. The
homotopy invariance of the index is shown directly in
Section~\ref{Homotopy} and in Section~\ref{Dirac} projective Dirac
operators are defined and the usual local index formula is used to compute
the index in that case. Much as in the usual case this formula is extended
to general pseudodifferential operators in Section~\ref{Index}. Some
examples in which the index is truly fractional are given in the final
section.

\section{Homomorphism bundles}\label{Hom}

Let $\cA$ be a finite-dimensional (star) Azumaya bundle over a compact
manifold $Z;$ see \cite{Mathai-Melrose-Singer1} for more details. By
definition $\cA$ is a complex vector bundle over $Z$ with fibres having
algebra structures and with local (algebra) trivializations as $N\times N$
matrix algebras. Since the automorphism group of the star-algebra of
$N\times N$ matrices is $\PU(N)$ (acting by conjugation) the bundle of all
such trivializations, $\cP,$ is a principal $\PU(N)$-bundle. From the
Azumaya perspective, the `trivial' case is where
$\cA\otimes\hom(A_1)=\hom(A_2)$ is `stably' the homomorphism bundle of a
vector bundle over $Z$ and this corresponds to the existence of a stable
lifting of $\cP$ to a $\UU(M)$-principal bundle.

A projective vector bundle, $E,$ over $Z$ can be defined (a different
initial approach is taken in \cite{Mathai-Melrose-Singer1}) as a
projection-valued section of $\cA\otimes\cK,$ for the algebra $\cK$ of
compact operators on some Hilbert space $\cH.$ Any projection in
$\cA_z\otimes\cK$ is of finite rank, so over a set in which $\cA$ is
trivial this yields a vector bundle. However, the phases of the transition
maps between trivializations are not determined, and cannot, in general, be
chosen to satisfy the cocycle condition, so in general these are not vector
bundles. The transpose Azumaya bundle $\cA^t$ is $\cA$ with multiplication
reversed and $\cA\otimes\cA^t$ is trivial as an Azumaya bundle, since it
has structure group acting through the adjoint representation,
$\PU(N)\longrightarrow \PU(N^2),$ which lifts canonically to a $\UU(N)$
action. For any two projective vector bundles $E$ and $F$ associated to
$\cA$ it follows that $\hom(E,F),$ since it is associated to
$\cA\otimes\cA^t,$ is a true vector bundle.

The lift of $\cA$ to an Azumaya bundle over $\cP$ is trivial, \ie is a
homomorphism bundle, and correspondingly the lift of a projective vector
bundle $E$ associated to $\cA$ to $\cP$ is a finite-dimensional subbundle,
$\tE\subset\bbC^N\otimes\cH,$ over $\cP$ which is equivariant for the
standard action of $\UU(N)$ on $\bbC^N,$ interpreted as covering the
$\PU(N)$ action on $\cP.$ Since the action of $\UU(N)$ on $\hom(\tE,\tF),$
for the lifts of any two projective vector bundles associated to $\cA,$ is
through conjugation we see that this is a bundle over $\cP$ invariant under
the $\PU(N)$ action and hence, again, we see that it descends to
$\hom(E,F),$ as a well-defined vector bundle on $Z.$ On the other hand the
`big' homomorphism bundle $\Hom(\tE,\tF)$ is only a projective vector
bundle over $Z^2;$ it is associated to the \emph{external} tensor product
$\cA\boxtimes\cA^t$ over $Z^2.$ Since at the diagonal it is a vector
bundle, reducing there to $\hom(E,F)$ it is reasonable to expect it be
represented by a vector bundle in a neighbourhood of the diagonal. For our
purposes it is vitally important that this extension be made in such a way
that the composition properties also extend.

For a given metric on $Z$ set 
\begin{equation}
N_\epsilon =\{(z,z')\in Z^2;d_g(z,z')<\epsilon \}.
\label{mms3.59}\end{equation}

The projective unitary group, $\PU(N),$ can be written as a quotient of the
group, $\SU(N),$ of unitary matrices of determinant one:
\begin{equation}
\bbZ_N\longrightarrow \SU(N)\longrightarrow \PU(N).
\label{mms3.72}\end{equation}
In the following result, which is the foundation of subsequent
developments, we use the discreteness of the fibres of \eqref{mms3.72}.

\begin{proposition}\label{mms3.4} Given two projective bundles, $E$ and
$F,$ associated to a fixed Azumaya bundle and $\epsilon >0$ sufficiently
small, the exterior homomorphism bundle $\Hom(\tE,\tF),$ descends from a
neighborhood of the diagonal in $\cP\times\cP$ to a vector bundle,
$\Hom^{\cA}(E,F),$ over $N_\epsilon$ extending $\hom(E,F).$ For any three such
bundles there is a natural associative composition law
\begin{multline}
\Hom^{\cA}_{(z'',z')}(F,G)\times \Hom^{\cA}_{(z,z'')}(E,F)\ni
(a,a')\longmapsto a\circ a'\in\Hom^{\cA}_{(z,z')}(E,G),\\ (z,z''),(z'',z')\in
N_{\epsilon/2}
\label{mms3.5}\end{multline}
which is consistent with the composition over the diagonal.
\end{proposition}

\begin{remark}\label{mms3.83} Applying this result to the projective vector
bundle given by $\Id_{\cA}\otimes\pi_1,$ where $\pi_1$ is the projection
onto the first basis element of $\cH,$ gives a bundle, which we denote
$\ecA,$ which extends $\cA$ from the diagonal to some neighborhood
$N_\epsilon$ and which has the composition property as in \eqref{mms3.5} 
\begin{equation}
\ecA_{(z'',z)}\times \ecA_{(z,z'')}\longrightarrow \ecA_{(z,z')}.
\label{mms3.84}\end{equation}
We regard this as the natural extension of $\cA.$
\end{remark}

\begin{proof} Consider again the construction of $\hom(E,F),$ always for
two projective bundles associated to the same Azumaya bundle, $\cA.$ The
dual bundle $\tE'$ is associated to the adjoint Azumaya bundle; or as a
subbundle of $\bbC^N\otimes\cH$ over $\cP$ it is associated with the
adjoint action of $\SU(N)$ on $\bbC^N.$ The external tensor product over
the product of $\cP$ with itself, as a bundle over $Z^2,$
$\tF\boxtimes\tE'$ is therefore a subbundle of
$\hom(\bbC^N)\otimes\hom(\cH)$ over $\cP\times\cP$ equivariant for the
action of $\SU(N)\times\SU(N)$ over $\PU(N)\times\PU(N).$ Restricted to the
diagonal, $\cP\times\cP$ has the natural diagonal subbundle $\cP.$ The
restriction of $\tF\boxtimes\tE'$ to this submanifold has a $\PU(N)$
action, and so descends to the bundle $\hom(E,F)$ over $\Diag\equiv Z,$
since for $A\in\PU(N)$ we can take the same lift $\tA$ to $\SU(N)$ in each
factor and these different diagonal lifts lead to the same operator through
conjugation.

In a neighborhood, $N_\epsilon,$ of the diagonal there is a corresponding
`near diagonal' submanifold of $\cP\times\cP;$ for instance we can extend
$\cP$ over the diagonal to a subbundle $\widetilde{\cP}\subset\cP\times\cP$
by parallel transport normal to the diagonal for some connection on $\cA.$
Now, any two points of $\widetilde{\cP}$ in the same fibre over a point in
$N_\epsilon$ are related by the action of $(A',A)\in\PU(N)\times\PU(N)$
where $A'A^{-1}$ is in a fixed small neighborhood of the diagonal, only
depending on $\epsilon.$ It follows from the discreteness of the quotient
$\SU(N)\longrightarrow \PU(N)$ that, for $\epsilon >0$ sufficiently small,
on lifting $A$ to $\tA\in \SU(N)$ there is a unique neighboring lift,
$\tA',$ of $A'.$ The conjugation action of these lifts on $\tF\boxtimes\tE'$
is therefore independent of choices, so defining $\Hom^{\cA}(E,F)$ over
$N_\epsilon.$ This bundle certainly restricts to $\hom(E,F)$ over the
diagonal.

In fact this construction is independent of the precise choice of
$\widetilde{\cP}.$ Namely if $\Hom(\tE,\tF)\equiv \tF\boxtimes\tE'$ is
restricted to a sufficiently small open neighborhood, $N,$ of $\cP$ as the
diagonal in $\cP\times\cP,$ then the part of $\PU(N)\times\PU(N)$ acting on
the fibres of $N$ lifts to act linearly on $\Hom(\tE,\tF),$ so defining
$\Hom^{\cA}(E,F)$ as a bundle over the projection of $N$ into $Z^2.$ It follows
that this action is consistent with the composition of $\Hom(\tE,\tF)$ and
$\Hom(\tF,\tG)$ for any three projective bundles associated to $\cA.$ This
leads to the composition property \eqref{mms3.5}.
\end{proof}

As a bundle over $\cP,$ the projective bundle $\tE$ can be given an
$\SU(N)$ invariant connection. A choice of such connections on $\tE$ and
$\tF$ induces, as in the standard case, a connection on $\Hom(\tE,\tF)$
over $\cP\times\cP$ and hence a connection on $\Hom^{\cA}(E,F)$ over
$N_\epsilon.$

\begin{remark}\label{mms3.85} Since $\cH$ is a fixed Hilbert space we can
also identify $\cK$ as a trivial bundle over $Z^2.$ The construction about
then identifies $\Hom^{\cA}(E,F)$ as a subbundle of $\ecA\otimes\cK,$ as a
bundle over $N_\epsilon,$ with the composition \eqref{mms3.5} induced from
\eqref{mms3.84}.
\end{remark}

\begin{remark}\label{mms3.58} A particularly important case of an Azumaya
bundle is the Clifford bundle on any oriented even-dimensional manifold (in
the odd-dimensional case the complexified Clifford bundle is not quite an
Azumaya bundle but rather the direct sum of two). On a manifold of
dimension $2n$ this is locally isomorphic to the algebra of $2^n\times 2^n$
matrices. Letting $\cP_{\CCl}$ be the associated principal
$\PU(2^n)$-bundle of trivializations, we call the trivial bundle,
$\bbC^{2^n},$ over $\cP_{\CCl}$ the \emph{projective spin bundle}; the
relationship to the usual spin bundle is explained in
Section~\ref{Trivial}. Proposition~\ref{mms3.4} and Remark~\ref{mms3.83},
give an extension of the Clifford bundle to a bundle, $\eCCl,$ in
neighborhood of the diagonal as the `big' homomorphism bundle of this
projective spin bundle. The discussion below shows that this allows us to
define a projective spin Dirac operator even when no $\spin$ (or even
$\spinC)$ structure exits. However, it is an element of an algebra of
`differential operators' which does not have any natural action.
\end{remark}

\section{Twisting by a line bundle over $\cP$}\label{P-twisting}

In Proposition~\ref{mms3.4} we have described a canonical extension of a given
Azumaya algebra $\cA$ to a bundle $\ecA$ near the diagonal. In general this
canonical extension, the existence of which is based on the discreteness of
the cover of $\PU(N)$ by $\SU(N),$ is not unique as an extension with the
composition property \eqref{mms3.84}. Rather, it is based on the selection,
natural as it is, of the trivial bundle $\tE=\cP\times\bbC^N$ with its
natural $\SU(N)$ action, as generating $\cA$ through $\hom(\tE).$ In this
section we consider the possibility of other choices of bundle in place of
$\tE$ and hence other extensions of $\cA.$

To motivate this discussion, consider the case of an Azumaya bundle which
is trivial, in the sense that it is isomorphic to $\hom(W)$ for an Hermitian
vector bundle $W.$ The frame bundle $\cW$ of $W$ is a principal
$\UU(N)$ bundle to which $\hom(W)$ lifts to be the trivial bundle of
$N\times N$ matrices on which $\UU(N)$ acts through conjugation. Thus the
center acts trivially, so $\hom(W)$ can also be identified with the trivial
bundle of $N\times N$ matrices over $\cP=\cW/\UU(1).$ The circle bundle,
$\cL,$ over $\cP$ with total space $\cW$ has an induced $\SU(N)$ action and
$W$ can be identified with the bundle $L_W\otimes\bbC^N$ over $\cP,$ where
$L_W$ is the line bundle corresponding to $\cL.$ Abstracting this situation
we arrive at the corresponding notion for a general Azumaya bundle.

\begin{definition}\label{mms3.73} A \emph{representing bundle} for a star
Azumaya bundle $\cA$ is a vector bundle $\tV$ over $\cP$ equipped with an
action of $\SU(N)$ which is equivariant for the $\PU(N)$ action on $\cP$
with the center acting as scalars and with an isomorphism, as bundles of
algebras, of $\cA$ and $\hom(\tV)$ as a bundle over the base.
\end{definition}

\noindent When appropriate we consider the unique $\UU(N)$ action on $\tV$
for which the center also acts as scalars and such that the only
elements of the center acting trivially are elements of
$\bbZ_N\subset\SU(N).$ Note that this is consistent with the `trivial' case
discussed above.

We consider two such representing bundles, $\tV_1,$ $\tV_2,$ to be
equivalent if there is a bundle isomorphism between them which intertwines
the $\SU(N)$ actions and projects to intertwine the isomorphisms with
$\cA.$ To understand the non-equivalent representing bundles we study the
line bundles on $\cP.$ The fibres of $\cP$ are diffeomorphic to $\PU(N)$ so
all line bundles over $\cP$ have flat connections over the fibres.

\begin{proposition}\label{mms3.86} (\cf Kostant \cite{Kostant1} and
Brylinski \cite{Brylinski3}). The total space of any line bundle on $\cP$
admits a `fibre holonomy' action by $\SU(N)$ which is equivariant for the
$\PU(N)$ action on $\cP,$ is linear between the fibres and in which the
centre acts as the fibre holonomy; this canonical 
$\SU(N)$ action is unique up to
conjugation by a bundle isomorphism.
\end{proposition}
\noindent As for representing bundles we consider the unique $\UU(N)$
action on the line bundle for which the center acts also acts as scalars
and such that the only elements of the center acting trivially are elements
of $\bbZ_N\subset\SU(N).$

\begin{proof} Let $L$ be a given line bundle on $\cP;$ choose some
connection, $\nabla,$ on it with curvature $\omega\in\CI(\cP;\Lambda^2).$
Each fibre of $\cP$ is diffeomorphic to $\PU(N).$ Since
$H^2(\PU(N),\bbR)=\{0\}$ the restriction of $\omega$ to each fibre is
exact. Thus we can find a smooth $1$-form, $\alpha,$ on $\cP$ such that
$\omega -d\alpha$ vanishes on each fibre. It follows that the connection
$\nabla-\alpha$ has vanishing curvature on each fibre of $\cP.$

Now, the $\SU(N)$ action on $L$ is given by parallel transport with respect
to such a connection. For each smooth curve $c:[0,1]\longrightarrow\PU(N)$
with $c(0)=\Id$ and each point $l\in L_p$ consider the curve $c(t)p$ in the
fibre through $p\in\cP$ and let $s(l)\in L_{c(1)p}$ be obtained by parallel
transport along $L$ over $c(t)p.$ This certainly gives a smooth map $s(c)$
on the total space of $L$ which is linear on the fibres. Furthermore, since
the curvature on each fibre vanishes, $s(c)$ depends only on the homotopy
class of $c$ in $\PU(N)$ as a curve from $\Id$ to $g=c(1)\in\PU(N)$ and
composition of curves leads to the composite map. Thus in fact $s$ is an
action of the universal covering group, $\SU(N),$ of $\PU(N)$ on the total
space of $L$ as desired. The centre, $\bbZ_N,$ of $\SU(N)$ gives the fibre
holonomy essentially by definition.

Any two connections on $L$ which are fibre-flat differ by a $1$-form $\beta$
which is closed on each fibre. Again, since $H^1(\PU(N),\bbR)=\{0\},$ we
may choose $f\in\CI(\cP)$ such that $\beta -df$ vanishes on each
fibre. Parallel transport along curves in $\PU(N)$ as discussed above, for
the two connections, is then intertwined by the bundle isomorphism
$\exp(f).$ Thus the $\SU(N)$ action defined by parallel transport on the
fibres is well-defined up to bundle isomorphism.
\end{proof}

\begin{remark} An alternate proof of Proposition \ref{mms3.86}  uses 
Cheeger-Simons characters, and will be described here. 
As in the proof of  Proposition \ref{mms3.86},
\begin{enumerate}
\item \label{one}given a 
line bundle $L$ over $\cP$, we can always find a connection with curvature 
$F$ with the property that $F$ restricted to the fibers is trivial.
\item \label{two}$\SU(N)$-actions on $L$ that cover the $\PU(N)$ action on $\cP$
are obtained from characters of $\bbZ_N$ 
via the holonomy of flat connections on line bundles along $\PU(N)$.
\item \label{three}Consider the exact sequence 
\begin{equation}\label{eqn:exact}
0\to H^1(\cP, \bbR)/H^1(\cP, \bbZ) \to {\check H}^2(\cP) \to A^2(\cP)\to 0
\end{equation}
where $ {\check H}^2(\cP)$ denotes the  Cheeger-Simons characters of 
1-cycles on $\cP$
and 
\begin{multline*}
A^2(\cP)= \Big\{(c_1, F);\text{where }F\text{ is a closed 2-form on }
\cP\text{ representing }c_1,\\
\text{which is  in the image of }H^2(\cP, \bbZ)\text{ in }H^2(\cP,\bbR)\Big\},
\label{mms3.101}\end{multline*}
\cf page 25, middle formula of equation (3.3), in \cite{Hopkins-Singer}.
\end{enumerate}
 
Take any pair $(c_1(L), F)$ as in \eqref{one} above. Then there exists a
Cheeger-Simons character $\chi : Z_1(\cP) \to \bbR/\bbZ$, whose value on a
closed curve is the holonomy of some connection with curvature equal to
$F$. Now the exact sequence \eqref{eqn:exact} when restricted to any fiber
of $\cP$ reduces to,
$$
 0\to {\check H}^2(\PU(N)) \to H^2(\PU(N), \bbZ) \to 0.
$$
Therefore the Cheeger-Simons characters of 1-cycles on $\PU(N)$, 
are in one-to-one correspondence with line bundles on $\PU(N)$.
So we deduce that the map from line bundles on $\cP$ to $\SU(N)$-actions on 
$L$ that cover the $\PU(N)$ action on $\cP$, is simply given by the map 
$$ H^2(\cP, \bbZ) \ni c_1 \to r^*c_1 \in   H^2(\PU(N), \bbZ) $$
where $r$ is the restriction map to any fiber. 
\end{remark}

One consequence of Proposition~\ref{mms3.86} is that any line bundle on $\cP$ is
necessarily the $N$th root of a line bundle on the base.

\begin{lemma}\label{mms3.75} If $\tL$ is a line bundle on $\cP$ for a
given Azumaya bundle then $\tL^{\otimes N}$ descends to a line bundle $L$ over
the base $Z.$
\end{lemma}

\begin{proof} It follows from the equivariance that the center $\bbZ_N$ of
$\SU(N)$ acts on the fibre $\tL_q,$ $q\in\cP,$ at each point as
multiplication by $N$th roots of unity. Thus, in the induced action on
$\tL^{\otimes N}$ the center acts trivially, so $\tL^{\otimes N}$ has an
induced $\PU(N)$-action over $\cP$ and so descends to a bundle on the
base.
\end{proof}

\begin{lemma}\label{mms3.77} Any line bundle over $\cP$ has a connection
with curvature which is the lift of the form $\frac1N\pi^*\omega_L$ from
the base where $\omega_L$ is the curvature on the base of a connection on
$L=\tL^{\otimes N}.$
\end{lemma}

\begin{proof} Following the discussion in the proof of
Proposition~\ref{mms3.86} a given line bundle $L$ carries a connection with
curvature, $\omega,$ which is $\PU(N)$-invariant and vanishes on the
fibres. Taking a smooth lift of any vector field, $v,$ to a vector field
$v^*$ on $\cP$ the form $i_{v^*}\omega$ is well-defined independent of the
lift, and closed on each fibre. Since $H^1(\PU(N),\bbR)=\{0\},$ it follows
that the connection may be further modified, by a smooth 1-form which
vanishes on each fibre, so that the curvature is a basic and
$\PU(N)$-invariant form, \ie is the lift of a form from the base.
Computing in any local trivialization of $\cP$ gives the curvature in terms
of the induced connection on the $N$th power.
\end{proof}

\begin{remark} Given a line bundle $L$ over a compact manifold $M,$ the
problem of finding an $N$-th root of $L$ which is also a line bundle on $M$
can be approached as follows. One can take the $N$-th root of the
transition functions of $L$ with respect to a good cover $U_a$,
$g_{ab}^{1/N} : U_{ab} \to \UU(1).$ On triple overlaps, this gives a cocycle
\begin{equation}
t_{abc} =  g_{ab}^{1/N}  g_{bc}^{1/N}  g_{ca}^{1/N}  : U_{abc} \to \ZZ_N
\end{equation}
where $\ZZ_N = ker(s),$ $s$ being the second homomorphism in the short
exact sequence
\begin{equation}\label{ce1}
\ZZ_N \to \UU(1) \to \UU(1)
\end{equation}
which is given by $s(z) = z^N.$ Then the obstruction to the existence of an
$N$-th root for $L$ is given by the connecting homomorphism in the
corresponding long exact sequence in cohomology
\begin{equation}
\cdots \to H^1(M, \underline{ \UU(1)})\stackrel{\beta}{\to} H^2(Z,\ZZ_N)\to
H^2(M,\underline{\UU(1)})\to H^2(M,\underline{ \UU(1)}) \to\cdots 
\end{equation}
\ie the obstruction class is $\beta(L) \in H^2(M,\ZZ_N).$ In fact,
$\beta(L) = c_1(\cL)\; ({\rm mod}\; N).$ If $L$ has an $N$th root $L_0$
then all of the other $N$th roots of $L$ are of the form $L_0 \otimes R$,
where $R$ is a line bundle on $Z$ such that 
$R^{\otimes N} = 1.$ Hence the set of $N$th roots of $L$ is a
$\bbZ_N$-affine space with associated vector space $H^1(M,\bbZ_N).$

Applying this to the case $M=\cP$ it follows that all of the line bundles
which have $N$th powers a given bundle $L$ over the base are of the form
$\tL \otimes R$, where $R^{\otimes N} = 1,$ so form a $\bbZ_N$-affine
space with associated vector space $H^1(\cP, \bbZ_N).$

Recall that a principal $\PU(N)$ bundle $\pi:\cP\to M$ has an invariant,
$t(\cP) \in H^2(M,\bbZ_N),$ which measures the obstruction to lifting $\cP$
to a principal $\SU(N)$ bundle. This obstruction is obtained via the
connecting homomorphism of the exact sequence in cohomology associated to
the short exact sequence of sheaves of groups on $M,$
$$
1\to \bbZ_N \to \underline{\SU(N)} \to \underline{\PU(N)}\to 1,
$$
namely $t(\cP) = \delta(\cP),$ where $\delta:H^1(M,\underline{\PU(N)})
\to H^2(M, \bbZ_N)$ and is the $\text{mod}\; N$ analogue of the 2nd
Stieffel-Whitney class.  Then a theorem of Serre (\cf \cite{Grothendieck})
asserts that given any class $t\in H^2(M,\bbZ_N),$ there is a
principal $\PU(mN)$ bundle $\pi':\cQ\to M$ (for some $m\in \mathbb N$) such
that $t(\cQ) \in H^2(M, \bbZ_{m N})$ maps to $t \in H^2(M, \bbZ_N)$
under the standard inclusion of the coefficient groups. 
Such a bundle $\cQ$ is by no means
unique.

In particular, given any line bundle $L$ on $M,$ by the theorem of Serre
\cite{Grothendieck}, discussed above we know that there is a principal
$\PU(mN)$ bundle $\pi':\cQ\to M$ (for some $m\in \mathbb N$) such that
$t(\cQ) = \beta(L) \in H^2(M,\bbZ_N).$ Since $t(\cQ)$ and $\beta(L)$ are
characteristic classes,
$$
0=\pi'^*(t(\cQ))= t(\pi'^*(\cQ)) = \beta(\pi'^*(L))\in H^2(\cQ, \bbZ_N),
$$
that is, there is a line bundle $\tilde L$ on $\cQ$ which isomorphic to an
$N$th root of the lifted line bundle $\pi'^*(L)$ on $\cQ$, \ie ${\tilde
L}^{\otimes N} \cong \pi'^*(L).$

Note that the exact sequence 
$$
1\to \bbZ \to \bbZ  \to \bbZ_N\to 1,
$$
where the middle arrow is multiplication by $N$, determines 
the change of coefficients long 
exact sequence in cohomology, where ${\mathfrak d} : H^2(M,\bbZ_N)
\to H^3(M,\bbZ)$ is one of the connecting homomorphisms. Then
${\mathfrak d} (t(\cP)) \in H^3(M,\bbZ)$ is  equal to the Dixmier-Douady
invariant, which measures the obstruction to lifting $\cP$ to a principal 
$\UU(N)$ bundle over $M$. Of course, this is a less stringent requirement.

\end{remark}

The proof of Proposition~\ref{mms3.4} applies just as well to the line
bundle $\tL\boxtimes\tL^{-1}$ over $\cP\times\cP.$

\begin{lemma}\label{mms3.87} If $\tL$ is a line bundle over $\cP,$ the
bundle $\tL\boxtimes\tL^{-1}$ descends from a neighborhood of the diagonal
submanifold of $\cP\times\cP$ to a well-defined line bundle $\hL$ over a
neighborhood of the diagonal in $Z^2.$
\end{lemma}

In general a line bundle over $\cP$ with its $\SU(N)$ action, and the
corresponding $\UU(N)$ action, represents a `partial trivialization' of the
Azumaya bundle. If the $\bbZ_N$ action on the fibres of $\tL,$ arising from
the centre of $\SU(N),$ is injective then in fact the circle bundle
associated to $\tL$ is a lift of the principal $\PU(N)$-bundle, $\cP,$ to a
principal $\UU(N)$-bundle. If, at the other extreme, this $\bbZ_N$ action
is trivial then $\tL$ is simply the lift of a line bundle from the base.

\begin{proposition}\label{mms3.76} Any representing bundle for an Azumaya
bundle is equivalent to $\tL\otimes\tE$ with the induced $\SU(N)$ action,
where $\tE$ is the trivial bundle with standard $\SU(N)$ action and $\tL$
is a line bundle on $\cP$ with its fibre-holonomy $\SU(N)$ action.
\end{proposition}

\begin{proof} Let $\tV$ be a representing bundle for the Azumaya bundle
$\cA.$ By assumption $\tV$ is a bundle over $\cP.$ Let $\cF$ be the frame
bundle for $\tV.$ This is the principal $\UU(N)$-bundle with fibre at a
point $p\in\cP$ the space of trivializations of $\tV_p.$ Now, as part of
the data of a representing bundle, we are given an identification of $\cA$
with $\hom(\tV)$ as a bundle over the base. Since a point of $\cP$ is an
identification of the fibres of $\cA$ with $N\times N$ matrices, this data
picks out a $U(1)$ subbundle $\cL\subset\cF,$ consisting of the
isomorphisms between $\tV_p$ and $\bbC^N$ which realize this identification
at that point. Since the equivariant $\UU(N)$ action on $\tV$ has center
acting as scalars, $\cL$ has an induced equivariant $\UU(N)$ action coming
from the equivariant $\UU(N)$ action on $\tV$ and the standard $\UU(N)$
action on the trivial bundle. If we let $\tL$ be the line bundle over $\cP$
associated to $\cL$ then it has a $\UU(N)$ action and the restriction of
this to $\SU(N)$ must be the $\SU(N)$ action. The frame bundle of
$\tL^{-1}\otimes\tV$ has a natural $\UU(N)$-invariant section over $\cP,$
so $\tV$ is equivalent, as a representing bundle, to $\tL\otimes\tE$ where
$\tE$ is the standard, trivial, representing bundle.
\end{proof}

In view of this result we generalize projective bundles slightly by
allowing twisting by line bundles over $\cP.$

\begin{definition}\label{mms3.81} For any Azumaya bundle $\cA$ and
line bundle, $\tL,$ over $\cP$ an associated ($\tL-)$twisted projective
bundle is a subbundle of $(\tL\otimes\bbC^N)\otimes\cH$ which is invariant
under the tensor product $\SU(N)$ action, arising from the $\SU(N)$ action
on $\tL$ and the standard $\SU(N)$ action on $\bbC^N$ interpreted as
covering the $\PU(N)$ action on $\cP.$
\end{definition}

\begin{proposition}\label{mms3.78} Any choice of representing bundle,
$\tV\equiv\tL\otimes\tE,$ for an Azumaya bundle $\cA$ over $Z$ gives rise
to a vector bundle $\ecA_{\tL}$ which is defined in a neighborhood of the
diagonal of $Z^2,$ extends $\cA=\hom(V)$ from the diagonal, has the
composition property \eqref{mms3.84} and lifts canonically to $\Hom(\tV)$
over a neighborhood of the diagonal on $\cP\times\cP;$ with its composition
maps \eqref{mms3.84} it is isomorphic to $\ecA\otimes\hL$ where $\tL$ is
given by Proposition~\ref{mms3.76}.
\end{proposition}

\begin{proof} The proof of Proposition~\ref{mms3.4} may be used directly,
since no use is made of the fact that the $\SU(N)$ action there is the
standard one.
\end{proof}

\begin{remark}\label{mms3.80} The same argument also gives an extension
$\Hom^{\cA,\tL}(E,F)$ of $\hom(E,F)$ for any two $\tL$-twisted projective
bundles $\tE$ and $\tF$ for the same line bundle over $\cP.$
\end{remark}

\begin{remark}\label{mms3.92} Applying Proposition~\ref{mms3.78} to the
`trivial' case of an Azumaya bundle $\cA=\hom(W)$ for a vector bundle $W$
we recover recover $\ecA_W=\Hom^{\cA,\tL}=\Hom(W)$ in a neighborhood of the
diagonal. 
\end{remark}

\section{Trivialization and spin structures}\label{Trivial}

Corresponding to \eqref{mms3.72}, the principal $\PU(N)$-bundle $\cP$ of
local trivializations of the Azumaya bundle $\cA$ may have a lift to a
principal $\SU(N)$-bundle, or to a principal $\UU(N)$-bundle,
\begin{equation}
\xymatrix{
\SU(N)\ar[r]\ar[d]&\cP_{\SU(N)}\ar[d]\\
\PU(N)\ar[r]&\cP
}\Mor
\xymatrix{
\UU(N)\ar[r]\ar[d]&\cP_{\UU(N)}\ar[d]\\
\PU(N)\ar[r]&\cP.
}
\label{mms3.60}\end{equation}
In either case the Dixmier-Douady invariant of $\cA$ vanishes; conversely the
vanishing of the invariant implies the existence of such a lifting to a
$\UU(pN)$-bundle for $\cA\otimes\hom(G)$ for some bundle $G$ (and so with
rank a multiple of $N.)$

Since it is of primary concern below, consider the special case of the
Clifford bundle. Namely, a choice of metric on $Z$ defines the bundle of Clifford
algebras, with fibre at $z\in Z$ the (real or complexified) Clifford algebra
\begin{equation}
\begin{gathered}
\Cl_z(Z)=\left(\bigoplus_{k=0}^\infty (T^*_zZ)^k\right)/\langle \alpha
\otimes\beta +\beta \otimes\alpha -2(\alpha ,\beta )_g,\ \alpha ,\beta
\in T^*_zZ \rangle,\\
\CCl_z(Z)=\bbC\otimes\Cl_z(Z).
\end{gathered}
\label{mms3.35}\end{equation}
If $\dim Z=2n,$ this complexified algebra is isomorphic to the matrix algebra on
$\bbC^{2^n}.$ A local smooth choice of orthonormal basis over an open set
$\Omega \subset Z$ identifies $T^*\Omega$ with $\Omega \times\bbR^{2n}$ and
so identifies $\Cl(\Omega)$ with $\Omega \times\Cl(\bbR^{2n})$ as Azumaya
bundles. Choosing a fixed identification of $\CCl(2n)$ with the algebra of
complex $2^n\times 2^n$ matrices therefore gives a trivialization of
$\CCl(Z),$ as an Azumaya bundle, over $\Omega.$ As noted in
\cite{Mathai-Melrose-Singer1}, its Dixmier-Douady invariant is $W_3(Z).$

In particular the Clifford bundle is an associated bundle to the metric
coframe bundle, the principal $\SO(2n)$-bundle $\cF,$ where the action of
$\SO(2n)$ on the Euclidean Clifford algebra $\Cl(2n)$ is through the spin
group. Thus, the spin group may be identified within the Clifford algebra as 
\begin{equation}
\Spin(2n)=\{v_1v_2\cdots v_{2k}\in\Cl(2n);v_i\in \bbR^{2n},\ |v_i|=1\}.
\label{mms3.61}\end{equation}
The non-trivial double covering of $\SO(2n)$ comes through the
mapping of $v$ to the reflection $R(v)\in\operatorname{O}(2n)$ in the plane
orthogonal to $v$  
\begin{equation}
p:\Spin(2n)\ni a=v_1\cdots v_{2k}\longmapsto R(v_1)\cdots R(v_{2k})=R\in\SO(2n).
\label{mms3.62}\end{equation}
Thus $\cP$ may be identified with the bundle associated to $\cF$ by the
action of $\SO(2n)$ on $\CCl(2n)$ (or in the real case $\Cl(2n))$ where $R$
in \eqref{mms3.62} acts by conjugation by $a$  
\begin{equation}
\Cl(2n)\ni b\longmapsto aba^{-1}\in\Cl(2n).
\label{mms3.63}\end{equation}
We therefore have a map of principal bundles 
\begin{equation}
\cF\longrightarrow \cP.
\label{mms3.64}\end{equation}
Recall that the projective spin bundle on $\cP$ is just the bundle
associated to the natural action of $\Cl(n)$ on itself; it can therefore be
identified with the trivial bundle over $\cP$ with an equivariant $\SU(N)$
action, where $N=2^n.$

Now, a spin structure on $Z,$ corresponds to an extension, $\cF_{S},$ of the
coframe bundle to a $\Spin$ bundle, 
\begin{equation}
\xymatrix{\Spin(2n)\ar[r]\ar[d]^p&\cF_{S}\ar[d]\\
\SO(2n)\ar[r]&\cF.}
\label{mms3.65}\end{equation}
Since $\Spin(2n)\subset\SU(N),$ where $\SU(N)\subset\CCl(2n),$ this in turn
gives rise to a lift of $\cP$ to a principal $\SU(N)$ bundle: 
\begin{equation}
\xymatrix{\SU(N)\ar[rr]&&\cP_{\SU(N)}\ar[dd]\\
\Spin(2n)\ar[r]\ar[d]^p\ar[u]&\cF_{S}\ar[d]\ar@{^{.}.{>}}[ur]\\
\SO(2n)\ar[r]&\cF\ar[r]&\cP.}
\label{mms3.66}\end{equation}
Thus the projective bundle naturally associated to the Clifford bundle can
reasonably be called the projective spin bundle since a spin structure on
the manifold gives a lift of $\tE\otimes M,$ where $M$ is the $\bbZ_2$
bundle given by the $\spin$ structure, to the usual spin bundle.

As in the standard case, the Levi-Civita connection induces a natural,
$\SU(N)$ equivariant, connection on the projective spin bundle over $\cP.$
We use this below to define the projective spin Dirac operator; a choice of
spin structure, when there is one, identifies it with the spin Dirac operator.

Note that similar remarks apply to a $\spinC$ structure on the manifold
$Z.$ The model group 
\begin{multline}
\SpinC(2n)=\{cv_1v_2\cdots v_{2k}\in\Cl(2n);v_i\in \bbR^{2n},\ |v_i|=1,
c\in\bbC,\ |c|=1\}\\
=(\Spin\times\UU(1))/\pm,
\label{mms3.70}\end{multline}
is a central extension of $\SO(2n),$ 
\begin{equation}
\UU(1)\longrightarrow \SpinC(2n)\longrightarrow \SO(2n)
\label{mms3.71}\end{equation}
where the quotient map is consistent with the covering of $\SO(2n)$ by
$\Spin(2n).$ 

Thus a $\spinC$ structure is an extension of the coframe bundle to a
principal $\SpinC(2n)$-bundle;
\begin{equation}
\xymatrix{\UU(1)\ar[r]\ar[d]&L\\
\SpinC(2n)\ar[r]\ar[d]&\cF_{L}\ar[d]\\
\SO(2n)\ar[r]&\cF.}
\label{mms3.67}\end{equation}
where $\cF_L,$ the $\SpinC(2n)$ bundle, may be viewed as a circle bundle
over $\cF.$ Since
$\SpinC(2n)\hookrightarrow \UU(N)$ (but is not a subgroup of $\SU(N))$ this
gives a diagram similar to \eqref{mms3.66} but lifting to a principal
$\UU(N)$ bundle
\begin{equation}
\xymatrix{\UU(N)\ar[rr]&&\cP_{\UU(N)}\ar[dd]\\
\SpinC(2n)\ar[r]\ar[d]^p\ar[u]&\cF_{L}\ar[d]\ar@{^{.}.{>}}[ur]\\
\SO(2n)\ar[r]&\cF\ar[r]&\cP.}
\label{mms3.69}\end{equation}
In this case the $\spinC$ bundle over $Z$ is the lift of $S\otimes L$ from
$\cP$ to $\cP_{\UU(N)}.$

Note that the existence of a spin structure on $Z$ is equivalent to the
condition $w_2=0.$ The Clifford bundle is then the homomorphism bundle of the
spinor bundle, so the existence of a spin structure implies the vanishing
of the Dixmier-Douady invariant of the Clifford bundle (which is $W_3,$ the
Bockstein of $w_2);$ the vanishing of $W_3$ is precisely equivalent to the
existence of a $\spinC$ structure (without any necessity for stabilization).

In the general case, even when $W_3\not=0$ and there is no $\spinC$
structure, we shall show below that we can still introduce the notion of a
`projective $\spinC$ Dirac operator' starting from the following notion.

\begin{definition}\label{mms3.79} On any even-dimensional, oriented
manifold a projective $\spinC$ structure is a choice of representing
bundle, in the sense of Definition~\ref{mms3.73}, for the complexified Clifford
bundle.
\end{definition}

Thus, by Proposition~\ref{mms3.76} such a representing bundle is always
equivalent to, and hence can be replaced by, $\tS\otimes\tL$ where $\tL$ is
a line bundle over the bundle of trivializations $\cP_{\CCl}$ of the
Clifford bundle and $\tS$ is the projective spin bundle. As remarked above,
this is consistent with the standard case in which there is a $\spinC$
structure and $\tV$ then descends to a bundle on the base.

\begin{remark} Any line bundle over $\cP_{\CCl}$ is necessarily a square
root of a line bundle from the base. This follows by restricting the line
bundle to the frame bundle, as a subbundle of $\cP_{\CCl},$ and the $\UU(N),$
$N=2^n,$ action to $\Spin(2n)$ showing that the centre acts through the
subgroup $\ZZ_2\subset\ZZ_N.$
\end{remark}

\section{Smoothing operators}\label{Smooth}

For two vector bundles $E$ and $F$ the space of smoothing operators
$\Psi^{-\infty}(Z;E,F)$ between sections of $E$ and sections of $F$ may
be identified with the corresponding space of kernels
\begin{equation}
\Psi^{-\infty}(Z;E,F)=\CI(Z^2;\Hom(E,F)\otimes\pi_R^*\Omega )
\label{mms3.7}\end{equation}
where the section of the density bundle allows invariant integration. Thus,
such kernels define linear maps $\CI(Z;E)\longrightarrow \CI(Z;F)$
through
\begin{equation}
Au(z)=\int_ZA(z,z')u(z').
\label{mms3.8}\end{equation}
Operator composition induces a product
\begin{equation}
\begin{gathered}
\Psi^{-\infty}(Z;F,G)\circ \Psi^{-\infty}(Z;E,F)\subset
\Psi^{-\infty}(Z;E,G),\\
A\circ B(z,z'')=\int_ZA(z,z')B(z',z'')
\end{gathered}
\label{mms3.9}\end{equation}
using the composition law \eqref{mms3.5}. The right density factor in $A$
is used in \eqref{mms3.9} to carry out the integral invariantly.

Given the extensions in Proposition~\ref{mms3.4} and
Proposition~\ref{mms3.78} of the homomorphism bundles it is possible to
define the linear space of smoothing operators with kernels supported in
$N_\epsilon$ for any pair $E,$ $F$ of projective bundles (or twisted
projective bundles) associated to a fixed Azumaya bundle (and twisting) as
\begin{equation}
\Psi^{-\infty}_{\epsilon}(Z;E,F)=
\CIc(N_\epsilon;\Hom^{\cA,\tL}(E,F)\otimes\pi_R^*\Omega)
\label{mms3.11}\end{equation}
where in case $E$ and $F$ are projective bundles, without twisting, $\tL$
is trivial so is dropped from the notation. Note that the projective and
possibly twisted nature of $E$ and $F$ is implicit in the
notation. Although there is no action analogous to \eqref{mms3.8} the composition
law \eqref{mms3.5} allows \eqref{mms3.9} to be extended directly to define
\begin{equation}
\Psi^{-\infty}_{\epsilon/2}(Z;F,G)\circ
\Psi^{-\infty}_{\epsilon/2}(Z;E,F)
\subset \Psi^{-\infty}_{\epsilon}(Z;E,G)
\label{mms3.12}\end{equation}
in the case of three projective bundles associated to the fixed $\cA.$ For
sufficiently small supports this product is associative
\begin{multline*}
(A\circ B)\circ C=A\circ(B\circ C)\\
\Mif A\in\Psi^{-\infty}_{\epsilon/4}(Z;G,H),\
B\in\Psi^{-\infty}_{\epsilon/4}(Z;F,G),\
C\in\Psi^{-\infty}_{\epsilon/4}(Z;E,F).
\label{mms3.13}\end{multline*}

The trace functional extends naturally to these spaces
\begin{equation}
\Tr:\Psi^{-\infty}_{\epsilon}(Z;E)=\int_Z \tr A(z,z)
\label{mms3.14}\end{equation}
and vanishes on appropriate commutators
\begin{equation}
\Tr(AB-BA)=0\Mif A\in\Psi^{-\infty}_{\epsilon/2}(Z;F,E),\
B\in\Psi^{-\infty}_{\epsilon/2}(Z;E,F)
\label{mms3.15}\end{equation}
as follows from Fubini's theorem.

\section{Pseudodifferential operators}\label{Pseudo}

Just as the existence of the bundle $\Hom^{\cA,\tL}(E,F)$ over the neighborhood
$N_\epsilon$ of the diagonal allows smoothing operators to be defined,
it also allows arbitrary pseudodifferential operators, with kernels
supported in $N_\epsilon$ to be defined as the space of kernels
\begin{equation}
\Psi^{m}_{\epsilon}(Z;E,F)=I^m_{\text{c}}(N_\epsilon;\Diag)\otimes_
{\CIc(N_\epsilon)}\CIc(N_\epsilon ;\Hom^{\cA,\tL}(E,F)).
\label{mms3.16}\end{equation}
Here, $E,F$ are either $\tL$-twisted projective bundles associated to 
some Azumaya bundle $\cA.$ These are just the conormal sections of
$\Hom^{\cA,\tL}(E,F)\otimes\pi_R^*\Omega$ with support in $N_\epsilon.$
Notice that for any small $\delta<\epsilon,$
\begin{equation}
\Psi^{m}_{\delta}(Z;E,F)+\Psi^{-\infty}_{\epsilon}(Z;E,F)=
\Psi^{m}_{\epsilon}(Z;E,F).
\label{mms3.20}\end{equation}

The singularities of these kernels are unrestricted by the support
condition so there are the usual short exact sequences
\begin{equation}
\begin{gathered}
\Psi^{m-1}_{\epsilon}(Z;E,F)\longrightarrow
\Psi^{m}_{\epsilon}(Z;E,F)\overset{\sigma _m}\longrightarrow
\CI(S^*Z;\hom(E,F)\otimes N_m)\\
\Psi^{-\infty}_{\epsilon}(Z;E,F)\longrightarrow
\Psi^{m}_{\epsilon}(Z;E,F)\overset{\sigma}\longrightarrow
\rho^{-m}\CI(S^*Z;\hom(E,F))[[\rho ]].
\end{gathered}
\label{mms3.17}\end{equation}
In the first case $N_m$ is the bundle over $S^*Z$ of the smooth functions
on $T^*Z\setminus0$ which are homogeneous of degree $m;$ this sequence is
completely natural and independent of choices. In the second sequence
$\rho\in\CI(\com{T^*Z})$ is a defining function for the boundary and the
image space represents Taylor series at the boundary, with an overall
factor of $\rho ^{-m};$ this sequence depends on choices of a metric and
connection to give a quantization map.

The product for pseudodifferential operators extends by continuity (using
the larger spaces of symbols with bounds, rather than the classical symbols
implicitly used above) from the product for smoothing operators and leads
to an extension of \eqref{mms3.12}
\begin{equation}
\begin{gathered}
\Psi^{m}_{\epsilon/2}(Z;F,G)\circ \Psi^{m'}_{\epsilon/2}(Z;E,F)
\subset \Psi^{m+m'}_{\epsilon}(Z;E,G),\\
\sigma _{m+m'}(A\circ B)=\sigma _{m'}(A)\sigma _m(B).
\end{gathered}
\label{mms3.18}\end{equation}
The induced product on the image space given by the second short exact
sequence in \eqref{mms3.17} is a star product as usual.

The approximability of general pseudodifferential operators by smoothing
operators, in the weaker topology of symbols with bounds, also shows, as in
the standard case, that \eqref{mms3.15} extends to
\begin{equation}
\Tr([A,B])=0\Mif A\in\Psi^{m}_{\epsilon/2}(Z;F,E),\
B\in\Psi^{-\infty}_{\epsilon/2}(Z;E,F)
\label{mms3.25}\end{equation}
for any $m.$

Now, the standard \emph{symbolic} constructions of the theory of
pseudodifferential operators carry over directly since these are all
concerned with the diagonal singularity and the symbol map.

\begin{theorem}\label{mms3.19} For any two projective bundles associated
to the same Azumaya bundle (or twisted projective bundles associated to the
same Azumaya bundle and the same line bundle over $\cP$), if
$A\in\Psi^{m}_{\epsilon/2}(Z;E,F)$ is elliptic, in the sense that $\sigma
_m(A)$ is invertible (pointwise), then there exists
$B\in\Psi^{-m}_{\epsilon/2}(Z;F,E)$ such that
\begin{equation}
B\circ A=\Id-E_R,\ A\circ B=\Id-E_L,\
E_R\in\Psi^{-\infty}_{\epsilon}(Z;E),\
E_L\in\Psi^{-\infty}_{\epsilon}(Z;F)
\label{mms3.21}\end{equation}
and any two such choices $B',$ $B$ satisfy
$B'-B\in\Psi^{-\infty}_{\epsilon/2}(Z;E,F).$
\end{theorem}

\begin{proof} Now absolutely standard.
\end{proof}

If $B'$ and $B$ are two such parametrices it follows that
$B_t=(1-t)B'+tB,$
$t\in[0,1],$ is a smooth curve of parametrices. Furthermore
\begin{equation}
\frac{d}{dt}[A,B_t]=[A,(B-B')]
\label{mms3.24}\end{equation}
so, by \eqref{mms3.25}, it follows that for any two parametrices
\begin{equation}
\Tr([A,B'])=\Tr([A,B])
\label{mms3.26}\end{equation}
since $B'-B$ is smoothing.

\begin{definition}\label{mms3.22} For an elliptic pseudodifferential
operator $A\in\Psi^{m}_{\epsilon/2}(Z;E,F)$ acting between projective
bundles associated to a fixed Azumaya bundle, or more generally between
twisted projective bundles corresponding to the same twisting line bundle
over $\cP,$ we define
\begin{equation}
\ind_a(A)=\Tr(AB-\Id_F)-\Tr(BA-\Id_E)
\label{mms3.23}\end{equation}
for any parametrix as in Theorem~\ref{mms3.19}.
\end{definition}

\section{Homotopy invariance}\label{Homotopy}

\begin{proposition}\label{mms3.27} The index \eqref{mms3.23} is constant
on a 1-parameter family of elliptic operators.
\end{proposition}

\begin{remark}\label{mms3.28} Given the rationality proved in the next
section this follows easily. Here we use the homotopy invariance to
prove the rationality!
\end{remark}

\begin{proof} For a smooth family
$A_t\in\CI([0,1];\Psi^m_{\epsilon/2}(Z;E,F))$ of elliptic operators as
discussed above, it follows as in the standard case that there is a
smooth
family of parametrices,
$B_t\in\CI([0,1];\Psi^{-m}_{\epsilon/2}(Z;E,F)).$
Thus the index, defined by \eqref{mms3.23} is itself smooth, since
$[B_t,A_t]$ is a smooth family of smoothing operators. To prove directly
that this function is constant we use the residue trace of Wodzicki, see
\cite{Wodzicki7}, with an improvement to the definition due to Guillemin
\cite{Guillemin2} and also the trace defect formula from
\cite{Melrose-Nistor2}.

For a classical operator $A$ of integral order, $m,$ in the usual
calculus the residue trace is defined by `$\zeta$-regularization'
(following ideas of Seeley) using the entire family of complex powers of a
fixed positive (so self-adjoint) elliptic operator of order $1:$
\begin{equation}
\rTr(A)=\lim_{z\to0}z\Tr(AD^{z})
\label{mms3.29}\end{equation}
where $\Tr(AD^z)$ is known to be meromorphic with at most simple poles at
$z=-k-\dim Z+\{0,1,2,\dots\}.$ One of Guillemin's innovations was to show
that the same functional results by replacing $D^z$ by any entire family
$D(z)$ of pseudodifferential operators of complex order $z$ which is
elliptic and has $D(0)=\Id.$

One way to construct such a family, which is useful below, is to choose a
generalized Laplacian on the bundle in question, which is to say a second
order self-adjoint differential operator, $L,$ with symbol $|\xi|^2\Id,$
the metric length function, and to construct its heat kernel, $\exp(-tL).$
This is a well-defined (locally integrable) section of the homomorphism
bundle on $[0,\infty)_t\times Z^2$ which is singular only at
$\{t=0\}\times\Diag(Z)$ and vanishes with all derivatives at $t=0$ away
from the diagonal. If $L=D^2$ is strictly positive the heat kernel decays
exponentially as $t\to\infty$ and the complex powers of $L$ are given by
the Mellin transform
\begin{equation}
L^z=D^{2z}=\frac1{\Gamma (z)}\int_0^{\infty} t^{-z-1}\exp(-tL)dt,
\label{mms3.44}\end{equation}
where the integral converges for $\Re z<<0$ and extends meromorphically to
the whole of the complex plane. The fact that $D^0=\Id$ arises from the
residue of the integral at $z=0,$ so directly from the fact that
$\exp(-tL)=\Id$ at $t=0.$ It follows that if $\chi \in\CIc([0,\infty)\times
Z^2)$ and $\chi \equiv1,$ in the sense of Taylor series, at
$\{0\}\times\Diag$ then
\begin{equation}
D(2z)=\frac1{\Gamma (z)}\int_0^{\infty} t^{-z-1}H(t)dt,\quad
H(t)=\chi\exp(-tL)
\label{mms3.45}\end{equation}
is an entire family as required for Guillemin's argument, that is $D$ is
elliptic and $D(0)=\Id.$ Since the construction of the singularity of the heat
kernel for such a differential operator is completely symbolic (see for
instance Chapter 5 of \cite{MR96g:58180}), quite analogous to the
construction of a parametrix for an elliptic operator, it can be carried out
in precisely the same manner in the projective case, so giving a family of
the desired type via \eqref{mms3.45}.

Alternatively, for any projective vector bundle $E$, such a family can be
constructed using an explicit linear quantization map, with kernels
supported arbitrarily close to the diagonal
\begin{equation}
D_E(z)\in\Psi^z_{\epsilon/4}(Z;E).
\label{mms3.30}\end{equation}
Thus we may define the residue trace and prove its basic properties as
in the standard case; in particular it vanishes on all operators of
sufficiently low order. It is also a trace functional
\begin{equation}
\rTr([A,B])=0,\ A\in \Psi^m_{\epsilon /4}(Z;E,F),\
B\in\Psi^{m'}_{\epsilon/4}(Z;F,E).
\label{mms3.31}\end{equation}

The additional result from \cite{Melrose-Nistor2}, see also
\cite{fipomb},
that we use here concerns the regularized trace. This is defined to be
\begin{equation}
\bTr_{D_E}(A)=\lim_{z\to0}\left(\Tr(AD_E(z))-\frac1z\rTr(A)\right).
\label{mms3.32}\end{equation}
For general $A$ it does depend on the regularizing family, but for
smoothing operators it reduces to the trace. Therefore
\begin{equation}
\ind_a(A) = \bTr_D([A, B]) = \bTr_{D_F}(AB-\Id_F)  -  \bTr_{D_E}(BA-\Id_E) 
\end{equation}
for an elliptic operator $A\in \Psi^m_{\epsilon /4}(Z;E,F)$ and 
$B\in\Psi^{m'}_{\epsilon/4}(Z;F,E)$ a parametrix for $A$.
It is not a trace function but
rather the `trace defect' satisfies
\begin{equation}
\bTr_D([A,B])=\rTr(B\delta _DA)
\label{mms3.33}\end{equation}
where $\delta _D : \Psi^\bullet_{\epsilon }(Z;E,F)\to \Psi^\bullet_{\epsilon }(Z;E,F)$
is defined as  $\delta _D A = \frac{d}{dz} \big|_{z=0}D_F(z) A D_E(-z)$ as in \cite{starprod}. 
Then one computes,
\begin{equation}
\begin{aligned}
\ind_a(A) = & \lim_{z\to 0}\left(\Tr(ABD_F(z))  - \Tr(BAD_E(z))\right) - \bTr_{D_F}(\Id_F)  + \bTr_{D_E}(\Id_E) \\
= &  \lim_{z\to 0}\left(\Tr(BD_F(z)A)  - \Tr(BAD_E(z))\right) - \bTr_{D_F}(\Id_F)  + \bTr_{D_E}(\Id_E) \\
= &  \lim_{z\to 0}\left(\Tr(B(D_F(z)AD_E(-z)-A)D_E(z)) \right)   - \bTr_{D_F}(\Id_F)  + \bTr_{D_E}(\Id_E) \\
= & \rTr(a^{-1}\delta_D a) - \bTr_{D_F}(\Id_F)  + \bTr_{D_E}(\Id_E) 
\end{aligned} 
\label{mms3.104}\end{equation}
where $a$ is the image of $A$ in the full symbol space and we observe that 
$D_E(z) D_E(-z) = \Id_E + O(z^2)$.
Now $\delta _D$ also satisfies
\begin{equation}
\rTr(\delta _Da)=0\ \forall\ a,
\label{mms3.46}\end{equation}
and when $E=F$, it is a derivation acting on the full symbol algebra in
\eqref{mms3.17}.

From these formul\ae\ the homotopy invariance of the index in the
projective case follows. Namely
\begin{equation}
\begin{aligned}
\frac{d}{dt}\ind_a(A_t)
=&\rTr(a_t^{-1}\delta_D\dot a_t)+\rTr((\frac{d}{dt}a_t^{-1})\delta
_Da_t,)\\
=&-\rTr(\dot a_t\delta _Da_t^{-1})-\rTr(a_t^{-1}\dot a_t a_t^{-1}\delta
_Da_t)=0.
\end{aligned}
\label{mms3.34}\end{equation}
Here, $a_t$ is the image of $A_t$ in the full symbol space in which
the
image of $B_t$ is $a_t^{-1}$ and \eqref{mms3.46} has been used.
\end{proof}

\begin{remark}\label{mms3.53} A similar argument also proves the
multiplicativity of the index. Thus if $A_i$ for $i=1,2$ are two elliptic
projective operators with the image bundle of the first being the same as
the domain bundle of the second, they can be composed if their supports are
sufficiently small. Let $B_i$ be corresponding parametrices, again with
very small supports. Then $B_1B_2$ is a parametrix for $A_2A_1$ and the
index of the product is given by \eqref{mms3.33} in terms of the `full
symbols' $a_i$ of the $A_i$ 
\begin{multline}
\ind_a(A_2A_1)=\bTr_D([A_2A_1,B_1B_2])=\rTr(a_1^{-1}a_2^{-1}\delta _D(a_2a_1))\\
=\rTr(a_1^{-1}\delta _Da_1)+\rTr(a_2^{-1}\delta _Da_2)=\ind_a(A_1)+\ind_a(A_2).
\label{mms3.54}\end{multline}
\end{remark}

\begin{remark} Another consequence of the homotopy invariance of the index
is, as noted after the definition, that it is necessarily real. We do not
use the reality in the proof of the index formula below, from which it
follows that the index is rational, so we only sketch the argument.

First observe that there is an elliptic operator of any order, on any
projective bundle, of index $0.$ Namely $D(m),$ discussed above, has this
property, since it commutes with the regularizing family $D(z)$ in the
symbol algebra, so (when $E=F)$ the index vanishes from
\eqref{mms3.104}. Thus, using the multiplicativity, we need only consider
the case of operators of order $0.$

Next note that, using inner products on the projective bundles,
$\inda(P^*)=-\ind(P)$ for any elliptic operator $P.$ To see this, consider
the operator on the direct sum of the bundles 
\begin{equation}
\tilde{P}=\begin{pmatrix}
0&P^*\\
-P&0
\end{pmatrix}.
\label{mms3.105}\end{equation}
If $Q$ is a parametrix for $P$ then 
\begin{equation*}
\begin{pmatrix}
0&-Q\\
Q^*&0
\end{pmatrix}.
\label{mms3.106}\end{equation*}
is a parametrix for \eqref{mms3.105}. Inserting this into the definition of
the index it follows directly that the index of $\tilde{P}$ is
$\inda(P)+\inda(P^*).$

Now, $\tilde{P}$ can be imbedded in the elliptic family 
\begin{equation}
\begin{pmatrix}
\sin(\theta)\Id_E&\cos(\theta)P^*\\
-\cos(\theta)P&\sin(\theta)\Id_F
\end{pmatrix}.
\label{mms3.107}\end{equation}
From the homotopy invariance it follows that the index is zero, thus indeed
$\inda(P)=-\inda(P^*).$ However, simply taking the complex conjugate of
\eqref{mms3.23} it then follows that 
\begin{equation}
\overline{\inda(P)}=-\inda(P^*)=\inda(P)\text{ is real.}
\label{mms3.108}\end{equation}
\end{remark}

\section{Projective Dirac operators}\label{Dirac}

The space of differential operators `acting between' two projective
bundles
associated to the same Azumaya algebra is well defined, since these are
precisely the pseudodifferential operators with kernels with supports
contained in the diagonal; we denote by $\Diff k(Z;E,F)$ the space of
these
operators of order at most $k.$

Of particular interest is that, in this projective sense, there is a
`spin Dirac operator' on every oriented even-dimensional compact
manifold. As discussed in Section~\ref{Trivial} above, the projective
bundle associated to the spin representation is the projective spin bundle
of $Z,$ which we denote by $S;$ if $Z$ is oriented
it splits globally as the direct sum of two projective bundles $S^\pm.$
There are natural connections on $\CCl(Z)$ and $S^\pm$ arising from the
Levi-Civita connection on $T^*Z.$ As  discussed in Proposition~\ref{mms3.4},
the homomorphism bundle of $S,$ which can be identified with $\CCl(Z),$ has
an extension to $\eCCl(Z)$ in a neighborhood of the diagonal, and this
extended bundle also has an induced connection. The projective 
spin Dirac operator may
then be identified with the distribution
\begin{equation}
\eth=\cl\cdot\nabla_{L}(\kappa_{\Id}),\ \kappa_{\Id}=\delta (z-z')\Id_S.
\label{mms3.36}\end{equation}
Here $\kappa_{\Id}$ is the kernel of the identity operator in
$\Diff*(Z;S)$
and $\nabla_{L}$ is the connection restricted to the left variables with
$\cl$ the contraction given by the Clifford action of $T^*Z$ on the
left. As in the usual case, $\eth$ is elliptic and odd with respect to
the $\bbZ_2$ grading of $S$ and locally the choice of a spin structure
identifies this projective spin Dirac operator with the usual spin Dirac
operator.

More generally we can consider projective twists of this projective spin Dirac
operator. If $E$ is any unitary projective vector bundle over $Z,$
associated to an Azumaya bundle $\cA$ and equipped with a Hermitian
connection then $S\otimes E$ is a projective bundle associated to
$\CCl(Z)\otimes\cA$ and it is a Clifford module in the sense that
$$
\CCl(Z)\subset\hom(S\otimes E).
$$
The direct extension of \eqref{mms3.36}, using the tensor product
connection, gives an element $\eth_{E}\in\Diff1(Z;S\otimes E)$ which is
again $\bbZ_2$ graded. In the special case that $S\otimes E$ is a bundle
a related construction is given by M. Murray and M. Singer
\cite{Murray-Singer1}.

The relation between the index of twisted, projective spin Dirac operators
(or more generally, projective elliptic operators) and 
the distribution index of transversally elliptic operators,  will be discussed
in a subsequent paper.

\begin{theorem}\label{mms3.38} The positive part,
$\eth^+_E\in\Diff1(Z;S^+ \otimes E,S^-\otimes E)$ of the projective 
spin Dirac operator twisted by a unitary projective vector bundle 
$E$, has index
\begin{equation}
\ind_a(\eth^+_E)=\int_Z\Ahat(Z)\wedge\Ch_{\cA}(E)
\label{mms3.39}\end{equation}
where $\Ch_{\cA}:K^0(Z;\cA)\longrightarrow H^{\even}(Z;\bbQ)$ is the
Chern
character in twisted K-theory.
\end{theorem}

\begin{proof} The proof via the local index formula, see
\cite{Berline-Getzler-Vergne1} and also \cite{MR96g:58180}, carries over to
the present case. As discussed in section \ref{Homotopy}, the truncated 
heat kernel $H(t),$ formally representing $\exp(-t\eth_E^2),$
near $\Diag_Z \times \{t=0\}$, is well-defined as a smooth kernel on
$Z^2\times(0,\infty),$ with values in $\Hom^{\CCl\otimes\cA}(S\otimes
E)\otimes\Omega _R,$ 
modulo an
element of $\dCI(Z^2\times[0,\infty);\Hom^{\CCl\otimes\cA}(S\otimes
E)\otimes\Omega _R);$ that is vanishing to all orders at $t=0.$ Then we
claim that the analogue of the McKean-Singer formula holds,
\begin{equation} \ind_a(\eth_E^+)=\lim_{t\downarrow0}\STr(H(t))
\label{mms3.40}\end{equation} where $\STr$ is the supertrace, the
difference of the traces on $S^+\otimes E$ and $S^-\otimes E.$ 
The local index formula, as a
result of rescaling, asserts the existence of this limit and its evaluation
\eqref{mms3.39}.

In the standard case the McKean-Singer formula \eqref{mms3.40}, for the
actual heat kernel, follows by comparison with the limit as $t\to\infty,$
which explicitly gives the index. Indeed then the function
$\STr(\exp(-t\eth_E^2))$ is constant in $t.$ In the present case the index
is defined directly through \eqref{mms3.49a} so the argument must be
modified. If $H^\pm(t)$ are the approximate heat kernels of
$\eth_E^-\eth_E^+$ and $\eth_E^-\eth_E^+$ respectively, then both approach
the identity as $t\downarrow0.$ Thus for smoothing operators $K^\pm$ on the
appropriate bundles, $H^\pm(t)K^\pm\longrightarrow K^\pm$ as smoothing
operators as $t\downarrow0.$ Thus, from the continuity of the trace on
smoothing operators, the index can be rewritten 
\begin{equation*}
\inda(\eth_E^+)=\lim_{t\downarrow0}\Tr\left((\eth_E^+B-\Id_F)H^-(t)\right)-
\lim_{t\downarrow0}\Tr\left((B\eth_E^+-\Id_E)H^+(t)\right)
\label{mms3.102}\end{equation*}
where $B$ is a parametrix for $\eth_E^+.$

For $t>0$ these approximate heat kernels are smoothing, so the terms can be
separated showing that 
\begin{multline}
\inda(\eth_E^+)=\lim_{t\downarrow0}\Tr\left(H^+(t)-H^-(t)\right)+
\lim_{t\downarrow0}\Tr_E\left(B(\eth^+H^+(t)-H^-(t)\eth_E^+)\right)\\
=\lim_{t\downarrow0}\STr(H(t)).
\label{mms3.103}\end{multline}
Here we use the fact that the difference $\eth^+H^+(t)-H^-(t)\eth_E^+$ is,
again by the (formal) uniqueness of solutions of the heat equation, a
smoothing operator which vanishes rapidly as $t\downarrow0.$ This term
therefore makes no contribution to the index and we recover \eqref{mms3.40}
and hence the local index formula for projective Dirac operators.
\end{proof}

Let $\cP$ be the principal $\PU(N)$ bundle associated to $\cA$ and $\cP'$
be the principal $\PU(N')$ bundle associated to $\CCl(Z)$, \cf
Section~\ref{Trivial}. Twisting by $\cP$ and $\cP'$-twisting line bundles
$\tL$ and $\tL'$ respectively, does not affect the local discussion, only
the final formula. Thus if $\tS_{\tL'}=\tS\otimes\tL'$ is a projective
$\spinC$ bundle in the sense of Definition~\ref{mms3.79} and
$\tE_{\tL}=\tE\otimes\tL$ is a projective vector bundle we may define the
twisted projective $\spinC$ Dirac operator on it by choice of an
$\SU(N')$-invariant connection on the twisting bundle $\tL'$ and
$\SU(N)$-invariant connection on the twisting bundle $\tL$, the Levi-Civita
connection on $\tS$ and an $\SU(N)$-invariant connection on $\tE.$ As usual
we think of these bundles, $S_{L'}, E_L$ as twisted projective bundles over
the manifold, although they are in fact bundles over $\cP'$ and $\cP$
respectively.  

\begin{theorem}\label{mms3.38a} The positive part,
$\eth^+_{L', E_L}\in\Diff1(Z; S^+_{L'}\otimes E_L, S^-_{L'}\otimes E_L)$ of the
projective $\spinC$ Dirac operator corresponding to a general projective
$\spinC$ structure and twisted by a unitary projective vector bundle $E$
has index
\begin{equation}
\ind_a(\eth^+_{L', E_L})=\int_Z\Ahat(Z)\wedge\exp(\frac12 c_1 ({L'}))
\wedge\Ch_{\cA}(E)\wedge\exp(\frac1N c_1(L))
\label{mms3.39a}\end{equation}
where $\Ch_{\cA}:K^0(Z;\cA)\longrightarrow H^{\even}(Z;\bbQ)$ is the
Chern character in twisted K-theory, $c_1(L)$ is the first Chern class of
$L,$  the $N$th power of the line bundle $\tilde L$ over $\cP$
and $c_1({L'})$ is the  first Chern class 
of $L',$ the square of the line bundle $\tilde L'$ over $\cP'.$
\end{theorem}

\section{Index formula}\label{Index}

\begin{theorem}\label{mms3.96} Given an Azumaya bundle, $\cA,$ over an even
dimensional compact manifold $Z$, the analytic index defines a map 
\begin{equation}
\inda:\Kc^0(T^*Z;\pi^*\cA)\longrightarrow \bbQ
\label{mms3.97}\end{equation}
where $\ind_a(A)=\ind_a(\sigma (A))$ for elliptic elements of
$\Psi_\epsilon (Z;E,F)$ for projective vector bundles associated to $\cA$
and  
\begin{equation}
\inda(b)= \int_{T^*Z}\Td(T^*Z)\wedge\Ch_{\cA}(b),\quad \forall\
b\in\Kc(T^*Z;\pi^*(\cA)). 
\label{mms3.98}\end{equation}
\end{theorem}

\begin{proof} It has been shown above that $\inda(A),$ for elliptic
elements of $\Psi_\epsilon (Z;E,F)$ is additive, homotopy invariant and
multiplicative on composition. Thus it does descend to a map as in
\eqref{mms3.97}, just as in the standard case, but with possibly real
values. As such a real-valued additive map on the twisted K-space
$\Kc(T^*Z;\pi^*(\cA)),$ $\inda$ must factor through the Chern character,
since it is an isomorphism over $\bbR$ (or $\bbQ).$ Thus  
\begin{equation}
\inda(b)=\widetilde{\inda}(\Ch_{\cA})(b)),\
\widetilde{\inda}:H^{\even}_{\text{c}}(T^*Z;\bbQ)\longrightarrow \bbR
\label{mms3.99}\end{equation}
being a well-defined map. However we may construct such elliptic projective
pseudodifferential operators by twisting the signature operator by a
projective vector bundle associated to $\cA.$ For these \eqref{mms3.39}
gives the index. From the Thom isomorphism in cohomology, we know that
these elements generate $H^{\even}_{\text{c}}(T^*Z;\bbQ)$ so suffice to
compute the map $\widetilde{\inda}$ in \eqref{mms3.99}. Thus it suffices to
show that the Riemann-Roch formula \eqref{mms3.39} is consistent with
\eqref{mms3.98}, but this follows from the standard case of the index
formula and linearity.
\end{proof}

In the non-oriented case we can pass to the oriented cover and deduce the
same formula. Similarly if we consider pseudodifferential operators acting
between $\tL$ twisted projective vector bundles corresponding to a line
bundle $\tL$ over the bundle of trivializations of an Azumaya bundle $\cA,$
and with $N$th power $L$ over the base, we arrive at the analogous twisted
formula generalizing \eqref{mms3.98} and \eqref{mms3.39a}
\begin{equation}
\inda(Q)=\int_{T^*Z}\Td(T^*Z)\wedge\Ch_{\cA}(\sigma
(Q))\wedge\exp(\frac1Nc_1(L)),\ \forall\ Q\in\Psi_{\epsilon}(Z;E,F)\text{
  elliptic.}
\label{mms3.100}\end{equation}

In the odd-dimensional case we may use suspension to reduce to the
even-dimensional case and again arrive at \eqref{mms3.100}. Namely take the
exterior tensor product with an untwisted operator of index one on the
circle. To do this it is necessary to generalize the discussion in
Section~\ref{Hom} to \ref{Homotopy} to such `product type' operators,
including the homotopy invariance, enough to show that this exterior tensor
product can be deformed, through elliptic operators in the product sense,
to a true (projective) elliptic pseudodifferential operator. This is
essentially a smooth analogue of arguments already present in
\cite{Atiyah-Singer3} and we forgo the details, since geometrically the
even dimensional case is the more interesting one.

\section{Fractions and the index formula}\label{Fractions}

On an oriented even-dimensional manifold, the vanishing of $W_3$ is
equivalent to the existence of a $\spinC$ structure  (\cf
\cite{Lawson-Michelsohn1}); in particular this follows if the manifold is
almost complex. In the almost complex case there is no $\spin$ structure
unless the canonical bundle has a square root. Nevertheless, there is
always a projective spin Dirac operator and Theorem \ref{mms3.38} applied
in this case case gives the usual formula
$$
\ind_a(\eth^+)=\int_Z\Ahat (Z).
$$
We recall some well known examples of oriented but non-spin manifolds
where $\int_Z\Ahat (Z)$ is a fraction, justifying the title of the paper.
The simplest is $Z = \bbC P^2$, in which case $\int_Z\Ahat (Z)=-\frac{1}{8}.$

Also in the almost complex case with Hermitian metric, we have the $\spinC$
Dirac operator 
\begin{equation}
\dbar+\dbar^*:\Lambda ^{0,\even}Z\longrightarrow \Lambda ^{0,\odd}Z.
\label{mms3.82}\end{equation}
Its index is $\int_Z\Ahat(Z)e^{\frac12 c_1}$ where $c_1=c_1(Z)$ is the
Chern class of the canonical line bundle. The integral is the formula for
the top term in the Todd polynomial written in terms of $\Ahat$ and $c_1.$

An amusing corollary of Theorem~\ref{mms3.38a} is that we can now interpret
the integral as the index of the projective Dirac operator coupled to a
line bundle which is a square root of the canonical bundle. Previously this
interpretation was only possible when $Z$ was itself spin, when this square
root bundle exists as an ordinary line bundle on $Z.$

Another important class of examples is the following. Let $V^{2n}(2d+1)$ be
hypersurfaces in $\bbC {\rm P}^{2n+1}.$ That is, in the homogeneous coordinates
$[Z_0, \ldots, Z_{2n+1}]$ for $\bbC {\rm P}^{2n+1}$,
$$
\begin{array}{lcl}
V^{2n}(2d+1) &=& \Big\{ [Z_0, \ldots, Z_{2n+1}]  \in \bbC P^{2n+1}: 
P(Z_0, \ldots, Z_{2n+1})=0, \\[+7pt]
& & \nabla P(Z_0, \ldots, Z_{2n+1}) \ne 0,
(Z_0, \ldots, Z_{2n+1})\ne 0 \Big\}
\end{array}
$$
where $P(Z_0, \ldots, Z_{2n+1})$ is a homogeneous polynomial of degree
$2d+1$. Then it is known that $V^{2n}(2d+1) $ is not a spin manifold, and
that
$$
\int_{V^{2n}(2d+1) }\Ahat (V^{2n}(2d+1) ) = \frac{2^{-2n} 
(2d+1)}{(2n+1)!} \prod_{k=1}^n ((2d+1)^2
- (2k)^2).
$$
It is straightforward to see that for $d\ge n$, the right hand side is
equal to a non-zero fraction that is not an integer.

Note that ${\mathbb C}{\rm P}^2$ has positive scalar curvature and the
Bochner-Lichnerowicz formula holds for the projective operator $\eth^2,$
yet $\ind_a(\eth) = \int_Z{\widehat A}(Z) = -\frac{1}{8}\ne 0$! The usual
argument, by contradiction, to the vanishing of the index, and hence
$\Ahat$ genus, is not applicable since in the twisted case there is no
notion of global section of the projective spinor bundle and therefore no
way to construct harmonic spinors.

As we have observed before, $Z$ has no $\spinC$ structure if $W_3(Z)\not=0.$
Nevertheless the projective Dirac operator exists and can have a nonzero
 index. We thank M.J. Hopkins for examples of $Z$ with both $W_3(Z) \ne 0$,
and $\int_Z\Ahat(Z) \not\in \bbZ$. Here is one of his examples. 
Let ${\rm S}^2\hookrightarrow \bbC {\rm P}^4$ be an embedding of degree $2.$ 
In homogeneous coordinates we can take the embedding to be 
$(x,y)\mapsto (x^2,y^2,xy,0,0)$. To do surgery on the embedded 
${\rm S}^2$, we need to verify that its complex normal bundle $N$ is trivial as a 
six dimensional real vector bundle $N_{\mathbb R}.$ It is not hard to show that 
$c_1(N) \in H^2(\bbC {\rm P}^4,\mathbb Z)\cong\mathbb Z$ is equal to $-4.$ The obstruction 
to $N_{\mathbb R}$ being isomorphic to ${\rm S}^2 \times {\mathbb R}^6$ is 
$w_2(N_{\mathbb R}).$ But $w_2(N_{\mathbb R})=c_1(N) \,{\rm mod}\, 2=-4\, {\rm mod}\, 2=0.$

We can now perform the surgery. A tubular neighborhood of the embedded ${\rm S}^2$ is 
${\rm S}^2 \times {\rm Disc}^6$ with boundary ${\rm S}^2 \times {\rm S}^5.$ 
Replace the tube by ${\rm Disc}^3 \times {\rm S}^5$ 
gluing its boundary ${\rm S}^2 \times {\rm S}^5$ to the tube boundary. 
We obtain a manifold $Z$ that is oriented cobordant to 
$\bbC {\rm P}^4$.
Hence $Z$ is not a spin manifold, i.e. $w_2(Z) \ne 0$. 
The surgery makes $H^2(Z,\mathbb Z)=0.$ 
Hence $W_3(Z) \ne 0$ from the usual long exact sequence, 
$$
\ldots \to H^2(Z,\bbZ)\to  H^2(Z,\bbZ) \to H^2(Z,{\bbZ}_2) \to H^3(Z,\bbZ) \to \ldots, 
$$
where the first arrow is multiplication by 2. Moreover,
$$
\int_Z\Ahat(Z)=\int_{\bbC {\rm P}^4}\Ahat(\bbC {\rm P}^4) =\frac3{128}.
$$ 


\begin{thebibliography}{10}

\bibitem{Atiyah-Singer3}
M.F. Atiyah and I.M. Singer, \emph{The index of elliptic operators, {III}},
  Ann. of Math. \textbf{87} (1968), 546--604, MR0236952, Zbl 0164.24301.

\bibitem{Berline-Getzler-Vergne1}
Nicole Berline, Ezra Getzler, and Mich{\`e}le Vergne, \emph{Heat kernels and
  {D}irac operators}, Springer-Verlag, Berlin, 1992 MR1215720, Zbl 0744.58001.

\bibitem{Brylinski3}
J-L. Brylinski, \emph{Loop spaces, characteristic classes and geometric
  quantization}, Progress in Mathematics, vol. 107, Birkh\"auser Boston, Inc.,
  Boston, MA, 1993, MR1197353, Zbl 0823.55002.


\bibitem{Grothendieck}
A. Grothendieck,  \emph{Le group de Brauer I, II, III},
 in `Dix Exposes sur la cohomologie des schemas', 
North Holland, Amsterdam, 1968, MR1608798, MR1608805, MR0244271, Zbl
0193.21503, Zbl 0198.25803, Zbl 0198.25901.


\bibitem{Guillemin2}
V.W. Guillemin, \emph{A new proof of {W}eyl's formula on the asymptotic
  distribution of eigenvalues}, Adv. Math. \textbf{55} (1985), 131--160,
MR0772612, Zbl 0559.58025.


\bibitem{Hopkins-Singer}
M. J. Hopkins and I.M. Singer, \emph{Quadratic functions in geometry, topology, 
and M-theory.} J. Differential Geom. \textbf{70} (2005), no. 3, 329--452,
MR2192936, 

\bibitem{Hormander3}
L.~H{\"o}rmander, \emph{The analysis of linear partial differential operators},
  vol.~3, Springer-Verlag, Berlin{,} Heidelberg{,} New York{,} Tokyo, 1985,
  MR1313500, Zbl 0601.35001.

\bibitem{Lawson-Michelsohn1}
H.~B.~Lawson Jr. and M.-L. Michelsohn, \emph{Spin geometry}, Princeton
  Mathematical Series, vol.~38, Princeton University Press, Princeton, NJ,
  1989, MR1031992, Zbl 0688.57001.
  
\bibitem{Kostant1}
Bertram Kostant, \emph{Quantization and unitary representations. I.
  prequantization.}, Lectures in modern analysis and applications, III, Lecture
  Notes in Math., vol. 170, Springer, Berlin, 1970, MR0294568, pp.~87--208,
MR0294568, Zbl 0223.53028.

\bibitem{Mathai-Melrose-Singer1}
V.~Mathai, R.B. Melrose, and I.M. Singer, \emph{The index of projective
  families of elliptic operators}, Geom. Topol. \textbf{9}, 341-373
(2005), MR2140985, Zbl pre02206530.

\bibitem{starprod}
Richard~B. Melrose, \emph{Star products and local line bundles}, 
\textbf{54}, No.5, 1581-1600 (2004), MR2127859, Zbl 1061.47064.


\bibitem{MR96g:58180}
Richard~B. Melrose, \emph{The {A}tiyah-{P}atodi-{S}inger index theorem}, A K
  Peters Ltd., Wellesley, MA, 1993, MR1348401, Zbl 0796.58050. 

\bibitem{Melrose-Nistor2}
Richard~B. Melrose and V.~Nistor, \emph{Homology of pseudodifferential operators on
  manifolds with corners {I}. {M}anifolds with boundary},
[{\tt{funct-an/9606005}}]

\bibitem{fipomb}
Richard~B. Melrose and {Fr\'ed\'eric} Rochon, \emph{Families index for
  pseudodifferential operators on manifolds with boundary},
Int. Math. Res. Not.  2004, \textbf{no. 22}, 1115--1141, MR2041651, Zbl
pre02207679.

\bibitem{Murray-Singer1}
Michael~K. Murray and Michael~A. Singer, \emph{Gerbes, clifford modules and the
  index theorem}, Ann. Global Anal. Geom. \textbf{26}, No.4, 355-367
(2004), MR2103405, Zbl pre02158449.

\bibitem{Wodzicki7}
M.~Wodzicki, \emph{Noncommutative residue. {I}. {F}undamentals}, {$K$}-theory,
  arithmetic and geometry (Moscow, 1984--1986), Lecture Notes in Math.,
  Springer, Berlin-New York, 1987, pp.~320--399, MR0923140, Zbl 0649.58033.

\end{thebibliography}
\def\cprime{$'$} \def\cdprime{$''$} \def\cprime{$'$} \def\cprime{$'$}
  \def\ocirc#1{\ifmmode\setbox0=\hbox{$#1$}\dimen0=\ht0 \advance\dimen0
  by1pt\rlap{\hbox to\wd0{\hss\raise\dimen0
  \hbox{\hskip.2em$\scriptscriptstyle\circ$}\hss}}#1\else {\accent"17 #1}\fi}
  \def\cprime{$'$} \def\ocirc#1{\ifmmode\setbox0=\hbox{$#1$}\dimen0=\ht0
  \advance\dimen0 by1pt\rlap{\hbox to\wd0{\hss\raise\dimen0
  \hbox{\hskip.2em$\scriptscriptstyle\circ$}\hss}}#1\else {\accent"17 #1}\fi}
  \def\cprime{$'$} \def\bud{$''$} \def\cprime{$'$} \def\cprime{$'$}
  \def\cprime{$'$} \def\cprime{$'$} \def\cprime{$'$} \def\cprime{$'$}
  \def\cprime{$'$} \def\cprime{$'$} \def\cprime{$'$} \def\cprime{$'$}
  \def\polhk#1{\setbox0=\hbox{#1}{\ooalign{\hidewidth
  \lower1.5ex\hbox{`}\hidewidth\crcr\unhbox0}}} \def\cprime{$'$}
  \def\cprime{$'$} \def\cprime{$'$} \def\cprime{$'$} \def\cprime{$'$}
  \def\cprime{$'$} \def\cprime{$'$} \def\cprime{$'$} \def\cprime{$'$}
  \def\cprime{$'$} \def\cprime{$'$} \def\cprime{$'$} \def\cprime{$'$}
  \def\cprime{$'$} \def\cprime{$'$} \def\cprime{$'$} \def\cprime{$'$}
  \def\cprime{$'$} \def\cprime{$'$} \def\cprime{$'$} \def\cprime{$'$}
  \def\cprime{$'$} \def\cprime{$'$} \def\cprime{$'$} \def\cprime{$'$}
  \def\cprime{$'$} \def\cprime{$'$}
\providecommand{\bysame}{\leavevmode\hbox to3em{\hrulefill}\thinspace}
\providecommand{\MR}{\relax\ifhmode\unskip\space\fi MR }
\providecommand{\MRhref}[2]{%
  \href{http://www.ams.org/mathscinet-getitem?mr=#1}{#2}
}
\providecommand{\href}[2]{#2}

\end{document}